\documentclass{amsart}
\usepackage{amscd,amssymb}
\usepackage{fullpage}
\usepackage{newlfont}
\usepackage{pb-diagram,lamsarrow,pb-lams}

\newcommand{\CK}        {{\mathcal{K}}}
\newcommand{\CM}        {{\mathcal{M}}}

\newcommand{\Aut}       {\operatorname{Aut}}
\newcommand{\End}       {\operatorname{End}}
\newcommand{\Ext}       {\operatorname{Ext}}
\newcommand{\Hom}       {\operatorname{Hom}}
\newcommand{\Map}       {\operatorname{Map}}

\newcommand{\Fp}        {{\mathbb{F}_p}}          
\newcommand{\Fq}        {{\mathbb{F}_q}}
\newcommand{\hFq}       {\widehat{\mathbb{F}}_q}
\newcommand{\Fqt}       {\mathbb{F}_q^\tm}
\newcommand{\Zp}        {{\mathbb{Z}_p}}          
\newcommand{\Zpt}       {{\mathbb{Z}_p^\tm}}
\newcommand{\Q}         {{\mathbb{Q}}}
\newcommand{\Qp}        {{\mathbb{Q}_p}}  
\newcommand{\QZp}       {{\mathbb{Q}_p/\mathbb{Z}_p}}              
\newcommand{\Z}         {{\mathbb{Z}}}
\newcommand{\C}         {{\mathbb{C}}}

\newcommand{\al}        {\alpha}
\newcommand{\bt}        {\beta} 
\newcommand{\gm}        {\gamma}
\newcommand{\dl}        {\delta}
\newcommand{\ep}        {\epsilon}
\newcommand{\tht}       {\theta}
\newcommand{\lm}        {\lambda}
\newcommand{\om}        {\omega}

\newcommand{\Gm}        {\Gamma}
\newcommand{\Lm}        {\Lambda}
\newcommand{\Sg}        {\Sigma}
\newcommand{\Om}        {\Omega}

\newcommand{\II}        {{\mathbb{I}}}
\newcommand{\Smash}     {\wedge}
\newcommand{\hSmash}    {\widehat{\wedge}}
\newcommand{\hI}        {\widehat{I}}
\newcommand{\hS}        {\widehat{S}}
\newcommand{\ot}        {\otimes}
\newcommand{\psb}[1]    {[\![#1]\!]}
\newcommand{\st}        {\;|\;}
\newcommand{\tm}        {\times}
\newcommand{\xra}       {\xrightarrow}
\newcommand{\xla}       {\xleftarrow}
\newcommand{\convto}    {\Longrightarrow}
\newcommand{\addet}     {\operatorname{addet}}
\newcommand{\hot}       {\widehat{\otimes}}
\newcommand{\cp}        {{\mathbb{C}P}} 
\newcommand{\cpi}       {{\mathbb{C}P^\infty}}
\newcommand{\ugm}       {\underline{\gamma}}
\newcommand{\bE}        {\overline{E}}

\newcommand{\GG}        {{\mathbb{G}}}
\newcommand{\OO}        {{\mathcal{O}}}
\newcommand{\bigWedge}  {\bigvee}
\newcommand{\iffa}      {\Leftrightarrow}
\newcommand{\sg}        {\sigma}
\newcommand{\spec}      {\operatorname{spec}}
\newcommand{\spf}       {\operatorname{spf}}
\newcommand{\res}       {\operatorname{res}}
\newcommand{\mxi}       {\mathfrak{m}}

\newcommand{\colim}  {\operatornamewithlimits{\underset{\longrightarrow}{lim}}}
\newcommand{\invlim} {\operatornamewithlimits{\underset{\longleftarrow}{lim}}}

\renewcommand{\:}{\colon}

\newtheorem{theorem}{Theorem}

\newtheorem{lemma}[theorem]{Lemma}
\newtheorem{proposition}[theorem]{Proposition}

\newtheorem{corollary}[theorem]{Corollary}
\theoremstyle{definition}
\newtheorem{remark}[theorem]{Remark}
\newtheorem{definition}[theorem]{Definition}

\begin{document}
\title{Gross-Hopkins duality}
\author{N.~P.~Strickland}
\date{\today}
\bibliographystyle{abbrv}

\maketitle

In~\cite{hogr:rap} Hopkins and Gross state a theorem revealing a
profound relationship between two different kinds of duality in stable
homotopy theory.  A proof of a related but weaker result is given
in~\cite{de:mmb}, and we understand that Sadofsky is preparing a proof
that works in general.  Here we present a proof that seems rather
different and complementary to Sadofsky's.  We thank I-Chiau Huang for
help with Proposition~\ref{prop-residue}, and John Greenlees for
helpful discussions. 

We first indicate the context of the Hopkins-Gross theorem.
Cohomological duality theorems have been studied in a number of
contexts; they typically say that 
\[ H^k(X^*)=H^{d-k}(X)^\vee \] 
for some class of objects $X$ with some notion of duality
$X\leftrightarrow X^*$ and some type of cohomology groups $H^k(X)$
with some notion of duality $A\leftrightarrow A^\vee$ and some integer
$d$.  For example, if $M$ is a compact smooth oriented manifold of
dimension $d$ we have a Poincar\'e duality isomorphism
\[ H^k(M;\Q)=\Hom(H^{d-k}(M;\Q),\Q) \] 
(so here we just have $M^*=M$).  For another example, let $S$ be a
smooth complex projective variety of dimension $d$, and let $\Om^d$ be
the sheaf of top-dimensional differential forms.  Then for any
coherent sheaf $F$ on $S$ we have a Serre duality isomorphism
\[ H^k(S;\Hom(F,\Om^d))=\Hom(H^{d-k}(S;F),\C). \]  
This can be seen as a special case of the Grothendieck duality theorem
for a proper morphism~\cite{ha:rd}, which is formulated in terms of
functors between derived categories.  There is a well-known analogy
between Boardman's stable homotopy category of spectra and the derived
category of sheaves over a scheme, and, inspired by this, Neeman has
used tools from homotopy theory to prove the main facts about
Grothendieck duality~\cite{ne:gdt}.  This leads one to expect that
there should be a kind of duality theorem in the stable homotopy
category itself.  However, experience suggests that such theorems
require finiteness conditions which are not satisfied in the category
of spectra.  The subcategory of $K(n)$-local
spectra~\cite{ra:lrc,host:mkl} has much better finiteness properties
and thus seems a better place to look for duality phenomena.  The
Hopkins-Gross theorem is a kind of analog of Serre duality in the
$K(n)$-local stable homotopy category.

We next explain the result in question in more detail; our formulation
will be compared with various other possible formulations in
Remark~\ref{rem-comparison}.  Fix a prime $p$ and an integer $n>0$ and
let $L_n$ be the Bousfield localisation functor~\cite{ra:lrc} with
respect to $E(n)$.  Let $MS$ be the fibre of the natural map
$L_nS\xra{}L_{n-1}S$, and let $\hI=IMS$ be its Brown-Comenetz dual,
which is characterised by the existence of a natural isomorphism
\[ [X,\hI] = \Hom(\pi_0(MX),\Q/\Z) = \Hom(\pi_0(MS\Smash X),\Q/\Z). \]
This lies in the $K(n)$-local category $\CK$, which is a symmetric
monoidal category with unit $\hS=L_{K(n)}S$ and smash product
$X\hSmash Y=L_{K(n)}(X\Smash Y)$.  The formal properties of the
category $\CK$ are studied in detail in~\cite{host:mkl}.  If $X\in\CK$
we write $\hI X=F(X,\hI)$ and $DX=F(X,\hS)$; these are the natural
analogues of the Brown-Comenetz dual and the Spanier-Whitehead dual of
$X$ in the $K(n)$-local context.  The following result is contained
in~\cite[Theorem 6]{hogr:rap}; it is convenient for us to state it
separately.
\begin{proposition}[Hopkins-Gross]\label{prop-invertible}
 The spectrum $\hI$ is invertible: there exists a spectrum
 $\hI^{-1}\in\CK$ with $\hI\hSmash\hI^{-1}=\hS$.  It follows that
 $\hI X=F(X,\hI)=\hI\hSmash DX$ for all $X\in\CK$.
\end{proposition}
\begin{proof}
 The original Hopkins-Gross proof has never appeared, but a proof of
 the first claim was given in~\cite[Theorem 10.2(e)]{host:mkl}.  The
 second follows purely formally from the first: the functor
 $\hI\hSmash(-)\:\CK\xra{}\CK$ is an equivalence, so we have
 \[ [Y,\hI X]=[X\hSmash Y,\hI]=[\hI^{-1}\hSmash X\hSmash Y,S]=
    [\hI^{-1}\hSmash Y,DX]=[Y,\hI\hSmash DX].
 \]
 Yoneda's lemma completes the argument.
\end{proof}

Our next problem is to study these dualities using a suitable
cohomology theory, which we call Morava $E$-theory.  It is a
$2$-periodic Landweber exact theory with coefficient ring 
\[ E_*=W\psb{u_1,\ldots,u_{n-1}}[u,u^{-1}] \]
where $|u_k|=0$ and $|u|=2$ and $W$ is the Witt ring of the finite
field $\Fq=\mathbb{F}_{p^n}$.  We make this a $BP_*$-algebra by the
map sending $v_k$ to $u^{p^k-1}u_k$, where $u_0=p$ and $u_n=1$ and
$u_k=0$ for $k>n$.  This gives rise to a Landweber exact formal group
law and thus a multiplicative cohomology theory.  We write $E$ for the
representing spectrum, which is $K(n)$-local, and is a wedge of
finitely many suspended copies of the spectrum $\widehat{E(n)}$.  If
$X\in\CK$ then we put $E^\vee_tX=\pi_t(E\hSmash X)$; this is the
natural version of $E$-homology to use in a $K(n)$-local context.  Our
goal is to describe $E^\vee_*\hI$.  The proposition above is actually
equivalent (by~\cite[Theorem 1.3]{homasa:cep}) to the statement that
$E^\vee_*\hI\simeq E_*$ as $E_*$-modules (up to suspension).  However,
this is not sufficient for the applications; we need pin down the
action of the group $\Gm$ of multiplicative automorphisms of $E$ (a
version of the Morava stabiliser group) as well as the $E_*$-module
structure.  This provides the input for various Adams-type spectral
sequences.  To explain the answer, we recall that there is a canonical
map $\det\:\Gm\xra{}\Zpt$ called the reduced determinant.  If $M$ is
an $E_0$-module with compatible action of $\Gm$ (in short, an
$E^0$-$\Gm$-module) then we write $M[\det]$ for the same $E^0$-module
with the $\Gm$-action twisted by $\det$.

\begin{theorem}[Hopkins-Gross]\label{thm-hopkins-gross}
 There is a natural $\Gm$-equivariant isomorphism
 \[ E_t^\vee\hI=E_{n-n^2+t}[\det]. \]
\end{theorem}

The proof of this result falls naturally into two parts.  One part is
the algebraic analysis of equivariant vector bundles on the Lubin-Tate
deformation space, which is outlined by Hopkins and Gross
in~\cite{hogr:rap}, with full details provided by the same authors
in~\cite{grho:evb}.  The other part is a topological argument to make
contact with Morava $E$-theory.  The theorem as stated implies a
description of $E_t^\vee\hI(Z)$ when $Z$ is a finite complex of type
$n$, and Devinatz~\cite{de:mmb} has given a proof (following a sketch
by Hopkins) that this description is valid when $p>(n^2+n+2)/2$.  

The rest of the paper will constitute our proof of the full theorem.

We first need some information about the the structure of the group
$\Gm$.  We will assume that the reader is familiar with the general
idea of the relationship between division algebras, automorphisms of
formal groups and cohomology operations.  We therefore give just
enough detail to pin down our group $\Gm$ among the various possible
versions of the Morava stabiliser group.  The original reference for
most of these ideas is~\cite{mo:nlc}, and~\cite{de:mcr} is a good
source for many technical points.

\begin{definition}
 Let $\phi$ be the unique automorphism of $W$ such that
 $\phi(a)=a^p\pmod{p}$ for all $a\in W$; note that $\phi^n=1$.  Let
 $D$ be the noncommutative ring obtained from $W$ by adjoining an
 element $s$ satisfying $sa=\phi(a)s$ (for $a\in W$) and $s^n=p$; thus
 $\Q\ot D$ is the central division algebra over $\Qp$ of rank $n^2$
 and invariant $1/n$.  There is a unique way to extend $\phi$ to an
 automorphism of $D$ with $\phi(s)=s$, and we still have $\phi^n=1$.
 This gives an action of the cyclic group
 $C=\langle\phi|\phi^n=1\rangle$ on $D$ and thus allows us to form the
 semidirect product $D^\tm\ltimes C$.
\end{definition}

\begin{proposition}\label{prop-gamma}
 There is a natural isomorphism $\Gm\simeq D^\tm\ltimes C$.
\end{proposition}
\begin{proof}
 Write $\Gm'=D^\tm\ltimes C$.  Let $F_0$ be the unique $p$-typical
 formal group law over $\Fq$ with $[p](x)=x^q$, and let $\GG_0$ be the
 formal group scheme over $S_0:=\spec(\Fq)$ associated to $F_0$.  If
 $\om\in\Fq=W/p$ then there is a unique lift $\hat{\om}\in W$ such
 that $\hat{\om}^q=\hat{\om}$, and there is also an endomorphism
 $\mu(\om)$ of $\GG_0$ given by $x\mapsto\om x$.  We also have an
 endomorphism $\sg$ of $\GG_0$ given by $x\mapsto x^p$.  It is
 well-known that there is a unique ring map
 $\tht\:D\xra{}\End(\GG_0)$ such that $\tht(s)=\sg$ and
 $\tht(\hat{\om})=\mu(\om)$ for all $\om\in\Fq$, and moreover that
 $\tht$ is an isomorphism.  On the other hand, we can let $C$ act on
 $S_0$ via the Frobenius automorphism $\phi\:S_0\xra{}S_0$.  As the
 coefficients of $F_0$ lie in $\Fp$, there is a natural identification
 $(\phi^k)^*\GG_0=\GG_0$.  Using this, it is easy to identify
 $\Gm'$ with the group of pairs $(\al_0,\bt_0)$, where
 $\bt_0\:S_0\xra{}S_0$ is an isomorphism of schemes and
 $\al_0\:\GG_0\xra{}\bt_0^*\GG_0$ is an isomorphism of formal groups
 over $S_0$. 

 Now put $S=\spf(E^0)$ and $\GG=\spf(E^0\cpi)$ so that $\GG$ is a
 formal group over $S$ and is the universal deformation of $\GG_0$ in
 the sense of Lubin and Tate~\cite{luta:fmo} (see also~\cite[Section
 6]{st:fsf} for an account in the present language.)  Let $\Gm''$
 denote the group of pairs $(\al,\bt)$ where $\bt\:S\xra{}S$ and
 $\al\:\GG\xra{}\bt^*\GG$.  As $S_0$ is the subscheme of $S$ defined
 by the unique maximal ideal in $E^0=\OO_S$ we see that
 $\bt(S_0)=S_0$.  As $\GG_0=S_0\tm_S\GG$, we also see that
 $\al(\GG_0)=\GG_0$, so we get a homomorphism $\Gm''\xra{}\Gm'$
 sending $(\al,\bt)$ to $(\al|_{S_0},\bt|_{\GG_0})$; deformation theory
 tells us that this is an isomorphism.  The general theory of
 Landweber-exact ring spectra~\cite{la:hpc} gives an isomorphism
 $\Gm''\xra{}\Aut(E)=\Gm$; see~\cite[Section 8.7]{st:fsf} for an
 account in the present language.
\end{proof}

For the next proposition, we note that the topological ring
$\Zp\psb{\Gm}$ has both a left and a right action of $\Gm$, and these
actions are continuous and commute with each other.  We can use the
left action to define continuous cohomology groups
$H^*(\Gm;\Zp\psb{\Gm})$, and the right action gives an action of $\Gm$
on these cohomology groups.
\begin{proposition}\label{prop-top-cohomology}
 We have $H^{n^2}(\Gm;\Zp\psb{\Gm})=\Zp$, and the natural action of
 $\Gm$ on this module is trivial.  Moreover, we have
 $H^k(\Gm;\Zp\psb{\Gm})=0$ for $k\neq n^2$.
\end{proposition}
\begin{proof}
 First, we note that $\Gm$ is a $p$-adic analytic group of dimension
 $n^2$ over $\Zp$.  Duality phenomena in the cohomology of profinite
 groups have been studied for a long time~\cite{se:cg,la:gap}, but the
 more recent paper~\cite{sywe:cpa} is the most convenient reference
 for the particular points that we need.  Write $U=1+pD<D^\tm<\Gm$,
 which is a torsion-free open subgroup of finite index in $\Gm$.  It
 follows from~\cite[Corollary 5.1.6]{sywe:cpa} that $U$ is a
 Poincar\'e duality group of dimension $n^2$ in the sense used in that
 paper, which means precisely that $U$ has cohomological dimension
 $n^2$ and $H^{n^2}(U;\Zp\psb{U})=\Zp$ and the other cohomology groups
 are trival.  As $U$ has finite index in $\Gm$, Shapiro's lemma gives
 an isomorphism $H^*(\Gm;\Zp\psb{\Gm})=H^*(U;\Zp\psb{U})$; this proves
 the proposition except for the fact that $\Gm$ acts trivially.  There
 is a unique way to let $\Gm$ act on $D$ such that the subgroup
 $D^\tm$ acts by conjugation and the subgroup $C$ acts via $\phi$.  On
 the other hand, $\Gm$ acts on $U$ by conjugation and thus on the
 $\Qp$-Lie algebra $L$ of $U$.  There is an evident $\Gm$-equivariant
 isomorphism $L=\Q\ot D$.  Using the results of \cite[Section
 5]{sywe:cpa} we get a $\Gm$-equivariant isomorphism
 \[ \Q\ot H^{n^2}(\Gm;\Zp\psb{\Gm})=\Q\ot H^{n^2}(U;\Zp\psb{U})=
     \Lm^{n^2}L=\Q\ot\Lm^{n^2}D.
 \]
 (This is just a tiny extension of an argument of Lazard, which could
 be applied directly if $\Gm$ were torsion-free.)  We write
 $\addet(\gm)$ for the determinant of the action of $\gm\in\Gm$ on
 $D$; it is now enough to check that $\addet=1$.

 Suppose $a\in D^\tm$, so $a$ acts on $D$ by $x\mapsto axa^{-1}$. Let
 $K$ be the subfield of $\Q\ot D$ generated over $\Qp$ by $a$.  Put
 $d=\dim_{\Qp}K$, so that $\Q\ot D\simeq K^{n^2/d}$ as left
 $K$-modules.  Using this, we see that the determinant of left
 multiplication by $a$ is just $N_{K/\Qp}(a)^{n^2/d}$.  By a similar
 argument, the determinant of right multiplication by $a^{-1}$ is
 $N_{K/\Qp}(a)^{-n^2/d}$.  It follows that the conjugation map has
 determinant one, so that $\addet(D^\tm)=1$.  Moreover, the action of
 $\phi\in C$ on $\Q\ot D$ is the same as the action of $s\in\Q\ot D$
 by conjugation, which has determinant one by the same argument.
\end{proof}

We next discuss some useful generalities about $E$-module spectra.  As
in~\cite[Appendix A]{host:mkl}, we say that an $E^0$-module $M$ is
\emph{pro-free} if it is the completion at $I_n$ of a free module.  
\begin{lemma}
 Let $M$ be a $K(n)$-local $E$-module spectrum.  Then the following
 are equivalent:
 \begin{itemize}
  \item[(a)] $\pi_1M=0$ and $\pi_0M$ is pro-free as an $E^0$-module.
  \item[(b)] $M$ is a coproduct in $\CK$ of copies of $E$. 
  \item[(c)] $M$ is a retract of a product of copies of $E$.
 \end{itemize}
 Moreover, if $\CM$ is the category of $E$-modules for which these
 conditions hold and $\CM_0$ is the category of pro-free $E$-modules
 then the functor $\pi_0\:\CM\xra{}\CM_0$ is an equivalence.
\end{lemma}
\begin{proof}
 Suppose~(b) holds, so $M=\bigWedge_\al E$ say.  Here the coproduct is
 the $K(n)$-localisation of the ordinary wedge, which means that
 $M=\invlim_J\bigWedge_\al E/J$, where $J$ runs over a suitable family
 of open ideals in $E^0$.  In this second expression, it makes no
 difference whether the wedge is taken in $\CK$ or in the category of
 all spectra, so $\pi_*\bigWedge_\al E/J=\bigoplus_\al E_*/J$, and it
 follows easily that $\pi_0M$ is the completion of $\bigoplus_\al E_0$
 and $\pi_1M=0$.  Thus (b)$\Rightarrow$(a).

 Conversely, if~(a) holds, choose a topological basis $\{e_\al\}$ for
 $\pi_0M$ and use it to construct a map $f\:\bigWedge_\al E\xra{}M$ of
 $E$-modules in the usual way.  Then $\pi_0\bigWedge_\al E$ is the
 completion of $\bigoplus_\al E_0$ and so $\pi_*(f)$ is an
 isomorphism, so~(b) holds.  Using this we find that for any
 $K(n)$-local $E$-module $N$, the group of $E$-module maps $M\xra{}N$
 is just $\prod_\al\pi_0N=\Hom_{E_0}(\pi_0M,\pi_0N)$.  If we write
 $\CM$ for the category of $E$-modules satisfying~(a) and~(b), it is
 now clear that $\pi_0\:\CM\xra{}\CM_0$ is an equivalence.
 
 We know from~\cite[Appendix A]{host:mkl} that $\CM_0$ is closed under
 arbitrary products and retracts, and that any pro-free $E^0$-module
 is a retract of a product of copies of $E^0$.  Given this, we can
 easily deduce that (a)$\iffa$(c).
\end{proof}
\begin{definition}
 We say that an $E$-module spectrum $M$ is \emph{pro-free} if it
 satisfies the conditions of the lemma.
\end{definition}
\begin{corollary}\label{cor-inv-lim}
 If $M$ is pro-free and $\{X_\al\}$ is the diagram of small spectra
 over $X$ then $[X,M]=\invlim_\al[X_\al,M]$.
\end{corollary}
\begin{proof}
 If $M=E$ then this holds by a well-known compactness argument based
 on the fact that $E^0Y$ is finite for all small $Y$.  As any $M$ can
 be written as a retract of a product of copies of $E$, the claim
 follows in general.
\end{proof}

\begin{lemma}\label{lem-smash-free}
 Let $M$ and $N$ be pro-free $E$-modules, and use the $E$-module
 structure on $M$ to make $M\hSmash N$ an $E$-module.  Then
 $M\hSmash N$ is pro-free.
\end{lemma}
\begin{proof}
 It is easy to reduce to the case $M=N=E$.  By the argument
 of~\cite[Proposition 8.4(f)]{host:mkl} it suffices to check that
 $(E/I_n)_*E$ is concentrated in even degrees.  We know by standard
 calculations that $(E/I_n)_*BP=E_*[t_k\st k>0]/I_n$ (with
 $|t_k|=2(p^k-1)$) and that $(E/I_n)_*E=(E/I_n)_*BP\ot_{BP_*}E_*$ by
 Landweber exactness, and the claim follows. 
\end{proof}

\begin{lemma}\label{lem-C-pro-free}
 Let $C(\Gm,E_0)$ be the ring of continuous functions from $\Gm$ (with
 its profinite topology) to $E_0$ (with its $I_n$-adic topology).
 Then $C(\Gm,E_0)$ is pro-free.
\end{lemma}
\begin{proof}
 Write $\hFq=\{a\in W\st a^q=a\}$; the reduction map $\hFq\xra{}\Fq$
 is well-known to be an isomorphism.  For any $d\in D^\tm$ there is a
 unique way to write $d=\sum_{k\geq 0}\tau_k(d)s^k$ with
 $\tau_k(d)\in\hFq$ and $\tau_0(d)\neq 0$.  This defines functions
 $\tau_k\:D^\tm\xra{}E_0$, and $\tau_k$ is constant on the cosets of
 $1+s^{k+1}D$ and thus is locally constant.  Lagrange interpolation
 shows that the evident map from $\Fq[\tau]/(\tau^q-\tau)$ to the ring
 of all functions $\Fq\xra{}\Fq$ is surjective, and thus an
 isomorphism by dimension count.  Our maps $\tau_i$ give a bijection 
 \[ (\tau_0,\ldots,\tau_k)\:D^\tm/(1+s^{k+1}D)
      \xra{}\Fqt\tm\mathbb{F}_q^k.
 \]
 Putting these facts together, we see that the ring of functions from
 $D^\tm/(1+s^{k+1}D)$ to $E_0$ is generated over $E_0$ by the
 functions $\tau_0,\ldots,\tau_k$ subject only to the relations
 $\tau_i^q=\tau_i$ and $\tau_0^{q-1}=1$.  The direct limit of these
 rings as $k$ tends to $\infty$ is the ring of all locally constant
 functions from $\Gm$ to $E_0$, which is thus isomorphic to
 \[ E_0[\tau_k\st k\geq 0]/(\tau_0^{q-1}-1,\tau_k^q-\tau_k). \]
 Now let $e_i\:C\xra{}E_0$ be the characteristic function of
 $\{\phi^i\}$ (for $i=0,\ldots,n-1$), so the ring of functions from
 $C$ to $E_0$ is just $E[e_0,\ldots,e_{n-1}]/(e_ie_j-\dl_{ij}e_i)$.

 Recall that $\Gm=D^\tm\ltimes C$; this can be identified with
 $D^\tm\tm C$ as a set, so the ring of locally constant functions from
 $\Gm$ to $E_0$ is just the tensor product of the rings for $D^\tm$
 and $C$, which we now see is a free $E_0$-module.  The ring of all
 continuous functions is the completion of the ring of locally
 constant functions, and thus is pro-free.
\end{proof}

We next want to identify $E^\vee_0E$ with $C(\Gm,E_0)$.  This is in
some sense well-known, but the details are difficult to extract from
the literature in a convenient form, so at the suggestion of the
referee we indicate a proof.  

It will first be helpful to have a more coherent view of the
topologies on our various algebraic objects.  Recall that a
$K(n)$-local spectrum $W$ is \emph{small} if it is a retract of the
$K(n)$-localisation of a finite spectrum of type $n$.  For any
$X,Y\in\CK$ we define a topology on $[X,Y]$ whose basic neighbourhoods
of $0$ are the kernels of maps $u^*\:[X,Y]\xra{}[W,Y]$ such that $W$
is small and $u\:W\xra{}X$.  We call this the ``natural topology'';
its formal properties are developed in~\cite[Section 11]{host:mkl}.
It is shown there that the natural topology on $E_0$ is the same as
the $I_n$-adic topology, which is easily seen to be the same as the
profinite topology.

\begin{lemma}\label{lem-Gamma-topology}
 The natural topology on $\Gm\subset[E,E]$ is also the same as the
 profinite topology.
\end{lemma}
\begin{proof}
 Consider a map $u\:W\xra{}E$ and the resulting map
 $u^*\:\Gm\xra{}E^0W$, which can also be thought of (using the evident
 action of $\Gm$ on $E^0W=[W,E]$) as the map $u\mapsto\gm^{-1}.u$.
 The set $E^0W$ is finite by an easy thick subcategory
 argument~\cite[Theorem 8.5]{host:mkl} so the stabiliser of $u$ has
 finite index in $\Gm$.  As the subgroup $1+pD<\Gm$ is a finitely
 topologically generated pro-$p$ group of finite index in $\Gm$, we
 see from~\cite[Theorem 1.17]{ddms:apg} that every finite index
 subgroup of $\Gm$ is open.  It now follows easily that the map
 $u^*\:\Gm\xra{}[W,E]$ is continuous when $[W,E]$ is given the
 discrete topology.  It follows in turn that if we give $\Gm$ the
 profinite topology and $E^0E$ the natural topology then the inclusion
 map is continuous.  We also know from~\cite[Proposition
 11.5]{host:mkl} that $E^0E$ is Hausdorff.  A continuous bijection
 from a compact space to a Hausdorff space is always a homeomorphism,
 and the lemma follows.
\end{proof}

\begin{theorem}\label{thm-coop}
 There is a natural isomorphism $E^\vee_0E=C(\Gm,E_0)$.
\end{theorem}
\begin{proof}
 Define a map $\phi\:\Gm\tm E^\vee_0E\xra{}E_0$ by
 \[ \psi(\gm,a) = 
   (S \xra{a} E\hSmash E \xra{1\hSmash\gm} E\hSmash E
    \xra{\text{mult}} E).
 \]
 Using~\cite[Propositions~11.1 and~11.3]{host:mkl} we see that $\psi$
 is continuous, so we have an adjoint map
 $\psi^\#\:E^\vee_0E\xra{}C(\Gm,E_0)$.  We claim that this is an
 isomorphism.  As both source and target are pro-free, it suffices to
 show that the reduction of $\psi^\#$ modulo $I_n$ is an isomorphism.
 Let $K$ be the representing spectrum for the functor
 $X\mapsto E^*/I_n\ot_{K(n)^*}K(n)^*X$; this is a wedeg of finitely
 many suspended copies of $K(n)$, and it can be made into an
 $E$-algebra spectrum with $K_*=E_*/I_n$.  It is not hard to see that
 \[ (E^\vee_0E)/I_n = (E/I_n)_0E = K_0E, \]
 and $C(\Gm,E_0)/I_n=C(\Gm,E_0/I_n)=C(\Gm,\Fq)$.  To analyse $K_0E$,
 let $x$ be the standard $p$-typical coordinate on $\GG$ and let $F$
 be the resulting formal group law over $E_0$ --- see~\cite[Appendix
 2]{ra:ccs} for details and useful formulae.  A standard Landweber
 exactness argument shows that $K_0E$ is the universal example of a
 ring $R$ equipped with maps $\Fq\xra{\al}R\xla{\bt}E_0$ and an
 isomorphism $f\:\bt_*F\xra{}\al_*F$ of formal group laws.  As
 $\al_*F$ has height $n$ we see that the same must be true of
 $\bt_*F$, so $\bt(I_n)=0$, so we can regard $\bt$ as a map
 $\Fq=E_0/I_n\xra{}R$.  The coefficients of $F$ modulo $I_n$ actually
 lie in $\Fp\leq\Fq$ and there is only one map $\Fp\xra{}R$ so
 $\al_*F=\bt_*F$; we just write $F$ for this formal group law.  Using
 the standard form for isomorphisms of $p$-typical FGL's we can write
 $f(x)=\sum^F_{k\geq 0}t_kx^{p^k}$, where $t_0$ is invertible because
 $f$ is an isomorphism.  (Readers may be more familiar with the graded
 case where one gets strict isomorphisms with $t_0=1$, but we are
 working with the degree zero part of two-periodic theories and $t_0$
 need not be $1$ in this context.)  As $f$ commutes with
 $[p]_F(x)=x^q$ we have
 $\sum^F_kt_kx^{p^{n+k}}=\sum^F_kt_k^qx^{p^{n+k}}$ and thus
 $t_k^q=t_k$.  In fact, this condition is sufficient for $f$ to be a
 homomorphism of FGL's (see~\cite[Appendix 2]{ra:ccs}, for example)
 and we deduce that
 \[ R=(\Fq\ot\Fq)[t_k\st k\geq 0]/(t_0^{q-1}-1,t_k^q-t_k)=
      C(D^\tm,\Fq\ot\Fq).
 \]
 We next claim that $\Fq\ot\Fq$ can be identified with the ring
 $F(C,\Fq)$ of functions from the Galois group
 $C=\langle\phi\st\phi^n=1\rangle$ to $\Fq$.  Indeed, we can define a
 map $\chi\:\Fq\ot\Fq\xra{}F(C,\Fq)$ by $\chi(a\ot b)(\sg)=a\sg(b)$.
 The $\Fq$-linear dual of this is the evident map
 $\Fq[C]\xra{}\End_{\Fp}(\Fq)$ which is injective by Dedekind's lemma
 on the independence of automorphisms, and this bijective by dimension
 count.  Thus $\chi$ is an isomorphism, and we obtain an isomorphism
 $R=C(D^\tm\ltimes C,\Fq)$.  After some comparison of definitions we
 see that this is the same as the map $\psi^\#\:K_0E\xra{}C(\Gm,\Fq)$,
 as required.
\end{proof}

\begin{definition}
 We write $J_k=E^{(k+1)}=E\hSmash\ldots\hSmash E$ (with $(k+1)$
 factors).  We can use the ring structure on $E$ to assemble these
 objects into a cosimplicial spectrum and thus a cochain complex of
 spectra.  We write $C_\Gm(\Gm^{k+1},E^0X)$ for the set of continuous
 $\Gm$-equivariant maps $\Gm^{k+1}\xra{}E^0X$ (where everything is
 topologised as in~\cite[Section 11]{host:mkl}).
\end{definition}

\begin{lemma}
 There is a natural isomorphism $[X,J_k]=C_\Gm(\Gm^{k+1},E^0X)$, which
 respects the evident cosimplicial structures.
\end{lemma}
\begin{proof}
 Given $\ugm=(\gm_0,\ldots,\gm_k)\in\Gm^{k+1}$, we define
 \[ \nu(\ugm) = (E^{(k+1)} \xra{\gm_0\Smash\ldots\Smash\gm_k}
                 E^{(k+1)} \xra{\text{mult}} E).
 \]
 We then define $\phi_X\:[X,J_k]\xra{}\Map(\Gm^{k+1},E^0X)$ by
 $\phi_X(a)(\ugm)=\nu(\ugm)\circ a\:X\xra{}E^0$.  Note that if
 $\dl\in\Gm$ then $\nu(\dl\ugm)=\dl\circ\nu(\ugm)\:J_k\xra{}E$; this
 means that the map $\phi_X(a)\:\Gm^{k+1}\xra{}E^0X$ is
 $\Gm$-equivariant.  The results of~\cite[Section 11]{host:mkl} show
 that $\phi_X(a)$ is continuous, so we can regard $\phi_X$ as a map
 $[X,J_k]\xra{}C_\Gm(\Gm^{k+1},E^0X)$.

 Now let $X$ be the dual of a generalised Moore
 spectrum~\cite{hosm:nshii,host:mkl} of type $S/I$ for some ideal
 $I=(u_0^{a_0},\ldots,u_{n-1}^{a_{n-1}})$ so that $E^0X=E_0/I$.  We
 know from Theorem~\ref{thm-coop} and Lemma~\ref{lem-C-pro-free} that
 $(E/I)_0E=C(\Gm,E_0/I)$ and that this is a free module over $E_0/I$
 so that
 \[ (E/I)_0E^{(k)}=C(\Gm,E_0/I)^{\ot k}=C(\Gm^k,E_0/I)=
     C_\Gm(\Gm^{k+1},E_0/I).
 \]
 After some comparison of definitions, we find that $\phi_X$ is an
 isomorphism when $X=D(S/I)$.  Moreover, the construction
 $M\mapsto C(\Gm^k,M)$ gives an exact functor from finite discrete
 Abelian groups to Abelian groups, so the construction
 $X\mapsto C(\Gm^k,E^0X)$ gives a cohomology theory on the
 category of small $K(n)$-local spectra.  By a thick subcategory
 argument, we deduce that $\phi_X$ is an isomorphism when $X$ is
 small.  Now let $X$ be a general $K(n)$-local spectrum, and let
 $\{X_\al\}$ be the diagram of small spectra over $X$.  As $J_k$ is a
 pro-free $E$-module we see from Corollary~\ref{cor-inv-lim} that
 $[X,J_k]=\invlim_\al[X_\al,J_k]$.  We also see that
 $C_\Gm(\Gm^{k+1},E^0X)=\invlim_\al C_\Gm(\Gm^{k+1},E^0X_\al)$, and it
 follows that $\phi_X$ is an isomorphism as claimed.
\end{proof}
\begin{corollary}\label{cor-ASS}
 There is a strongly convergent spectral sequence
 \[ E_1^{st}=[X,J_s]^t=C_\Gm(\Gm^{s+1},E^tX)\convto [X,\hS]^{t+s}. \]
\end{corollary}
\begin{proof}
 The axiomatic treatment of Adams resolutions discussed
 in~\cite{mi:ras} can be transferred to many other triangulated
 categories; this will certainly work for unital stable homotopy
 categories in the sense of~\cite{hopast:ash}, and thus for $\CK$.
 The spectra $J_s$ clearly form an $E$-Adams resolution of $\hS$, so
 we have a spectral sequence whose $E_1$ page is as described.  Let
 $i\:\bE\xra{}\hS$ be the fibre of the unit map $\eta\:\hS\xra{}E$; it
 is known that our spectral sequence is associated to the filtration
 of $\hS$ by the spectra $Y_s:=\bE^{(s)}$.  The map
 $\eta\hSmash 1\:E\xra{}E\hSmash E$ is split by the product map
 $\mu\:E\Smash E\xra{}E$, so $i\hSmash 1=0\:\bE\hSmash E\xra{}E$.  It
 follows that if $Z$ lies in the thick subcategory generated by $E$
 then the map $i^{(s)}\hSmash 1\:Y_s\hSmash Z\xra{}Z$ is zero for
 $s\gg 0$.  However, we know from~\cite[Theorem 8.9]{host:mkl} that
 $\hS$ lies in this thick subcategory, so $i^{(s)}=0$ for $s\gg 0$.
 Using this and the definition of our spectral sequence, we see that
 it converges strongly to $[X,\hS]^*$.
\end{proof}

\begin{proposition}\label{prop-DE}
 We have $[E,S^{n^2-t}]=E_t$ as $E_0$-modules, and thus
 $DE=\Sg^{-n^2}E$ as $E$-module spectra. Moreover, this equivalence
 respects the evident actions of $\Gm$.
\end{proposition}
\begin{proof}
 We use the spectral sequence of Corollary~\ref{cor-ASS}.  To analyse
 the $E_1$ term, we claim that $E^0E=\Zp\psb{\Gm}\hot_{\Zp}E_0$ as
 left $\Gm$-modules, where we let $\Gm$ act in the obvious way on
 $\Zp\psb{\Gm}$ and trivially on $E_0$.  To see this, define maps
 $\phi\:E_0\ot\Zp[\Gm]\xra{}E^0E$ and
 $\psi\:\Zp[\Gm]\ot E_0\xra{}E^0E$ by 
 \begin{align*}
  \phi(a\ot[\gm])  &= (E \xra{\gm} E \xra{\tm a} E) \\
  \psi([\gm]\ot a) &= (E \xra{\tm a} E \xra{\gm} E).
 \end{align*}
 Clearly $\gm'.\phi(a\ot[\gm])=\phi(\gm'(a)\ot[\gm'\gm])$ and
 $\gm'.\psi([\gm]\ot a)=\psi([\gm'\gm]\ot a)$ and
 $\psi([\gm]\ot a)=\phi(\gm(a)\ot[\gm])$ (because $\gm$ is a ring
 map).  By dualising Theorem~\ref{thm-coop} we see that $\phi$ extends
 to give an isomorphism $E_0\hot_{\Zp}\Zp\psb{\Gm}\xra{}E^0E$, and it
 follows from the above formulae that $\psi$ extends to give an
 isomorphism $\Zp\psb{\Gm}\hot_{\Zp}E_0\xra{}E^0E$ with the required
 equivariance.  Next, using the obvious description of
 $E_0=W\psb{u_1,\ldots,u_{n-1}}$ in terms of monomials, we see that
 $E_0$ is isomorphic as a topological $\Zp$-module to a product of
 copies of $\Zp$, say $E_0=\prod_\al\Zp$.  It follows that
 $E^0E=\prod_\al\Zp\psb{\Gm}$ as $\Gm$-modules, and using this that
 \[ C_\Gm(\Gm^{k+1},E^0E)=\prod_\al C_\Gm(\Gm^{k+1},\Zp\psb{\Gm})=
     C_\Gm(\Gm^{k+1},\Zp\psb{\Gm})\hot_{\Zp}E_0.
 \]
 
 These identifications can easily be transferred to nonzero degrees,
 and they respect the cosimplicial structure, so it now
 follows from Proposition~\ref{prop-top-cohomology} that in our
 spectral sequence we have
 \[ E_2^{s,t} = \begin{cases}
     E^t & \text{ if } s=n^2 \\
     0   & \text{ otherwise. }
    \end{cases}
 \]
 As the spectral sequence is strongly convergent, we have
 $[E,S]^r=E^{r-n^2}=E_{n^2-r}$ as claimed.  It is not hard to check
 that this is an isomorphism of $E_0$-modules, and it follows in the
 usual way that $DE=\Sg^{-n^2}E$ as $E$-module spectra.  One can see
 from the construction that this is compatible with the action of
 $\Gm$. 
\end{proof}

\begin{proposition}\label{prop-EvhI}
 There is a natural isomorphism $E^\vee_t\hI=\II(\pi_{n^2-t}ME)$. 
\end{proposition}
\begin{proof}
 First, we have
 \begin{align*}
  \II(ME_{t+n^2}\hI) &= [\Sg^{-n^2-t}E,\hI^2\hS] \\
    &= [\Sg^{-n^2-t}E,\hS] \\
    &= E_{-t}.
 \end{align*}
 We can deduce from the above by a thick subcategory argument that
 $K_*\hI$ is finite in each degree and thus that $\hI$ is dualisable.
 (We could also have quoted this from~\cite{host:mkl}; the proof given
 there is only a slight perturbation of what we've just done.)  This
 implies that $\hI X=F(X,\hI)=DX\hSmash\hI$ for all $X\in\CK$.  In
 particular, we have 
 \[ E\hSmash\hI=\Sg^{n^2}DE\hSmash\hI=\Sg^{n^2}\hI E, \]
 so that $E^\vee_t\hI=\pi_{t-n^2}\hI E=\II(\pi_{n^2-t}ME)$ as claimed.
\end{proof}

We next need to recall the calculation of $\pi_*ME$ and its connection
with local cohomology.  As usual, we define $E_*$-modules
$E_*/I_k^\infty$ and $u_k^{-1}E_*/I_k^\infty$ by the following
recursive procedure: we start with $E_*/I_0^\infty=E_*$, we define
$u_k^{-1}E_*/I_k^\infty$ by inverting the action of $u_k$ on
$E_*/I_k^\infty$, and we define $E_*/I_{k+1}^\infty$ to be the
cokernel of the evident inclusion
$E_*/I_k^\infty\xra{}u_k^{-1}E_*/I_k^\infty$.  It is also well-known
that there is a unique way to make $\Gm$ act on all these modules that
is compatible with the $E_*$-module structure and the short exact
sequences 
\[ E_*/I_k^\infty \xra{} u_k^{-1} E_*/I_k^\infty
   \xra{} E_*/I_{k+1}^\infty.
\]
(The point is that inverting $u_k$ is the same as inverting $v_k$, and
for each finitely generated submodule $M_*<E_*/I_k^\infty$ there is
some $N$ such that $v_k^{p^N}$ acts $\Gm$-equivariantly on $M_*$; this
determines the action on $u_k^{-1}E_*/I_k^\infty$.)  

One can define $E$-module spectra $E/I_k^\infty$ and
$u_k^{-1}E/I_k^\infty$ and cofibrations
\[ E/I_k^\infty \xra{} u_k^{-1}E/I_k^\infty \xra{} E/I_{k+1}^\infty \]
such that everything has the obvious effect in homotopy; this can
either be done by Bousfield localisation~\cite{ra:lrc} or by the
theory of modules over highly structured ring spectra~\cite[Section
6]{grma:cat}.  The former approach has the advantage that it is
manifestly $\Gm$-equivariant.  It is easy to see from the definitions
that $M(u_k^{-1}E/I_k^\infty)=0$ so that
$M(E/I_k^\infty)=\Sg^{-1}M(E/I_{k+1}^\infty)$, and also that
$M(E/I_n^\infty)=E/I_n^\infty$; it follows that
$ME=\Sg^{-n}E/I_n^\infty$.  All this is well-known and we record it
merely to fix the grading conventions.

We next consider some parallel facts in local cohomology.  Recall that
if $R$ is a Noetherian local ring with maximal ideal $\mxi$ and $M$ is
an $R$-module then the local cohomology groups $H^i_{\mxi}(M)$ may be
defined as $\colim_k\Ext^i_R(R/\mxi^k,M)$; see~\cite{brhe:cmr} (for
example) for an account of these groups.  Clearly $\Aut(R)$ acts
naturally on $H^*_{\mxi}(R)$, and if any element $u\in\mxi$ acts
invertibly on $M$ then $H^i_{\mxi}(M)=0$.  In particular, the element
$u_k\in I_n$ acts invertibly on $u_k^{-1}E_*/I_k^\infty$, so
$H_{I_n}^*(u_k^{-1}E_*/I_k^\infty)=0$.  Thus, the short exact sequence
\[ E_*/I_k^\infty \xra{} u_k^{-1} E_*/I_k^\infty
   \xra{} E_*/I_{k+1}^\infty
\]
gives rise to a $\Gm$-equivariant isomorphism 
\[ H^s_{I_n}(E_*/I_k^\infty)\simeq H^{s-1}_{I_n}(E_*/I_{k+1}^\infty).
\]
It is also standard that $H^0_{I_n}(E_*/I_n^\infty)=E_*/I_n^\infty$
and the other local cohomology groups of this module vanish.  We thus
end up with a $\Gm$-equivariant isomorphism
\[ H^n_{I_n}E_* = E_*/I_n^\infty, \]
and the other local cohomology groups of $E_*$ vanish.  (This can be
obtained more directly using an appropriate stable Koszul complex, but
$\Gm$ does not act on that complex and the whole question of
equivariance is rather subtle from that point of view.)

The next ingredient that we need is a certain ``residue map''
\[ \rho\: E_0/I_n^\infty\ot_{E_0}\Om^{n-1} \xra{}\QZp. \]
Here $\Om^{n-1}$ is the top exterior power of the module of K\"ahler
differentials for $E_0$ relative to $W$; it is freely generated over
$E_0$ by the element $\ep:=du_1\wedge\ldots\wedge du_{n-1}$.  The
definition of $\rho$ involves the map $\tau\:W\xra{}\Zp$ that sends
$a\in W$ to the trace of the map $x\mapsto ax$, considered as a
$\Zp$-linear endomorphism of $W$.  This induces a map
$W/p^\infty\xra{}\Z/p^\infty$ which we also call $\tau$.  We define
$\tau'\:W\xra{}\Hom(W,\Zp)$ by $\tau'(a)(b)=\tau(ab)$; it is a
standard fact of algebraic number theory that this is an isomorphism.

Given a multiindex $\al=(\al_1,\ldots,\al_{n-1})\in\Z^{n-1}$ we define
$u^\al=\prod_iu_i^{\al_i}$; we also write $\tht$ for the multiindex
$(1,\ldots,1)$.  The group $E_0/I_n^\infty\ot_{E_0}\Om^{n-1}$ is a
direct sum of copies of $W[1/p]/W=W\ot\QZp$ indexed by monomials
$u^{-\tht-\al}\ep$ for which $\al_i\geq 0$ for all $i$.  We put
\[ \rho(\sum_\al a_\al u^{-\tht-\al} \ep) = \tau(a_0). \]
\begin{proposition} \label{prop-residue}
 The map $\rho$ is $\Gm$-invariant.
\end{proposition}
\begin{proof}
 This is an instance of the well-known principle of invariance of
 residues.  Most of the formulae involved are very old, and were
 originally interpreted in a complex analytic context; starting in the
 1960's they were transferred into algebraic
 geometry~\cite{ha:rd,li:dsd} but only in the very recent
 paper~\cite{hu:pmz} do they appear in the particular technical
 context that we need.  First note that $W$ is the smallest closed
 subring of $E_0$ containing all the $(q-1)$'th roots of unity in
 $E_0$, so it is preserved by $\Aut(E_0)$ and in particular by $\Gm$.
 Section~5 of the cited paper gives a canonical map
 \[ \res\:H^{n-1}_{I_n}(E_0/p^\infty\ot_{E_0}\Om^{n-1})
      \xra{}W/p^\infty,
 \]
 and the short exact sequence $E_0\xra{}E_0[1/p]\xra{}E_0/p^\infty$
 gives an isomorphism
 \[ E_0/I_n^\infty\ot_{E_0}\Om^{n-1}=H^n_{I_n}(\Om^{n-1})=
     H^{n-1}_{I_n}(E_0/p^\infty\ot_{E_0}\Om^{n-1}),
 \]
 and $\tau$ gives a canonical map $W/p^\infty\xra{}\QZp$.  By putting
 all these together, we get a map
 $\rho'\:E_0/I_n^\infty\ot_{E_0}\Om^{n-1}\xra{}\QZp$.  All the
 constructions involved are functorial and thus $\Gm$-invariant.  By
 examining the formulae in~\cite[Section 5]{hu:pmz} we see that
 $\rho=\rho'$. 
\end{proof}

\begin{proposition}\label{prop-IIpiME}
 There is a natural isomorphism
 $\II(\pi_tME)=\Om^{n-1}\ot_{E_0}E_{-n-t}$.
\end{proposition}
\begin{proof}
 This is essentially the local duality theorem (see~\cite[Theorem
 5.9]{hu:pmz} for example).  We first translate the claim using the
 isomorphism $ME=\Sg^{-n}E/I_n^\infty$; it now says that
 $\II(E_t/I_n^\infty)=\Om^{n-1}\ot E_{-t}$.  Using the evident
 isomorphism $E_{-t}=\Hom_{E_0}(E_t,E_0)$ we can reduce to the case
 $t=0$.  In that case the construction above gives a $\Gm$-invariant
 element $\rho\in\II(E_0/I_n^\infty\ot\Om^{n-1})$ and thus a $\Gm$
 equivariant map $\phi\:\Om^{n-1}\xra{}\II(E_0/I_n^\infty)$, defined
 by $\phi(\lm)(a)=\rho(a\lm)$.  This satisfies
 \[ \phi\left(\sum_\bt b_\bt u^\bt\ep\right)
        \left(\sum_\al a_\al u^{-\tht-\al}\right) =
              \sum_\al \tau(a_\al b_\al).
 \]
 As the map $a\ot b\mapsto\tau(ab)$ is a perfect pairing, it is easy
 to conclude that $\phi$ is an isomorphism.
\end{proof}

We next recall the theorem of Hopkins and Gross, which identifies
$\Om^{n-1}$ in terms of certain modules whose $\Gm$-action is easier
to understand.  First, we write $\om=E_2$, considered as an
$E_0$-module with compatible $\Gm$-action.  To explain the notation,
note that $J:=\widetilde{E}^0\cpi$ is the group of formal functions on
$\GG$ that vanish at zero.  Thus $J/J^2$ is the cotangent space to
$\GG$ at zero, which is naturally isomorphic to the group of invariant
one-forms on $\GG$, which is conventionally denoted by $\om$.  On the
other hand $J/J^2$ is also naturally identified with
$\widetilde{E}^0\cp^1=\widetilde{E}^0S^2=E_2$, so the notation is
consistent.

\begin{theorem}[Gross-Hopkins]\label{thm-gross-hopkins}
 There is a natural $\Gm$-equivariant isomorphism
 $\Om^{n-1}\simeq\om^{\ot n}[\det]$.
\end{theorem}
\begin{proof}
 See~\cite[Corollary 3]{hogr:rap} (which relies heavily
 on~\cite{grho:evb}). 
\end{proof}

We can now give our proof of the main topological duality theorem.
\begin{proof}[Proof of Theorem~\ref{thm-hopkins-gross}]
 This follows from Propositions~\ref{prop-EvhI} and~\ref{prop-IIpiME}
 and Theorem~\ref{thm-gross-hopkins} after noting that
 $\om^{\ot n}=E_{2n}$.
\end{proof}

\begin{remark}\label{rem-comparison}
 We conclude by comparing our formulation of the theorem with various
 other possibilities considered by other authors.  Firstly, there are
 several alternative descriptions of the spectrum $\hI$: we have
 \[ \hI = I M S = L_{K(n)} I S = L_{K(n)} I L_{K(n)} S. \]
 To see this, let $X$ be a finite spectrum of type $n$, so that
 $L_{K(n)} X=MX=L_nX$, and smashing with $X$ commutes with all the
 functors under consideration.  The spectra listed above are all
 $K(n)$-local, and the canonical maps between them can be seen to
 become isomorphisms after smashing with $X$, and the claim follows.

 Next, there are several alternatives to our spectrum $E$, most
 importantly the spectrum $E'$ constructed using the Witt ring of
 $\overline{\mathbb{F}}_p$ rather than $\Fq$.  This is in some ways a
 more canonical choice, but it has the disadvantage that $\pi_0E'$ is
 not compact (although it is linearly compact) and $\Aut(E')$ is not a
 $p$-adic analytic group.  In any case, $E'$ is a pro-free $E$-module,
 so it is not hard to translate information between the two theories
 if desired.

 Finally, some other authors use the functor
 $\invlim_IE'_0(S/I\Smash X)$ in place of $E^\vee_0X$.  Hopkins writes
 this as $\CK(X)$, and calls it the Morava module of $X$.  One can
 show that
 \[ \CK(X)=\pi_0(L_{K(n)}(E'\Smash X))=E'_0\hot_{E_0}E^\vee_0X \]
 for large classes of spectra $X$, for example all $K(n)$-locally
 dualisable spectra, or spectra for which $K(n)_*X$ is concentrated in
 even degrees.  However, this equality can fail (for example when
 $X=\bigWedge_IS/I$) and the formal properties of the functor
 $E^\vee_0X$ seem preferable in general.
\end{remark}


\end{document}